\documentclass[11pt, leqno, a4paper]{amsart}
\usepackage{amsmath,amstext,amsthm,amsfonts}
\usepackage{amssymb,latexsym}
\usepackage[cp1252]{inputenc}
\usepackage[T1]{fontenc}
\usepackage{lmodern}
\usepackage{color}
\usepackage[english]{babel}

\theoremstyle{plain}
\newtheorem{theorem}{\textbf{Theorem}}[section]
\newtheorem{thm}[theorem]{\textbf{Theorem}}
\newtheorem{prop}[theorem]{\textbf{Proposition}}
\newtheorem{lem}[theorem]{\textbf{Lemma}}
\newtheorem{cor}[theorem]{Corollary}
\newtheorem*{thm*}{\textbf{Theorem}}

\theoremstyle{definition}
\newtheorem{def-d}{\textbf{Definition}}

\newtheorem{rem}[theorem]{\textbf{Remark}}

\numberwithin{equation}{section}

\def\boxit#1{\leavevmode\hbox{\vrule\vtop{\vbox{\kern.33333pt\hrule
    \kern1pt\hbox{\kern1pt\vbox{#1}\kern1pt}}\kern1pt\hrule}\vrule}}

\newcommand{\ie}{\emph{i.e. }}

\newcommand{\noi}{\noindent}

\newcommand{\noid}{\noindent $\diamond$~}
\newcommand{\R}{\mathbb{R}}

\newcommand{\C}{\mathbb{C}}
\newcommand{\D}{\mathbb{D}}

\newcommand{\Cb}{\mathbb{C}^{\bullet}}
\newcommand{\Z}{\mathbb{Z}}
\newcommand{\N}{\mathbb{N}}

\newcommand{\chih}{\widehat{\chi}}
\newcommand{\Sh}{\widehat{S}}
\newcommand{\Mh}{\widehat{M}}

\newcommand{\gh}{\hat{g}}

\newcommand{\4}{\frac{1}{4}}
\newcommand{\2}{\frac{1}{2}}
\newcommand{\lip}{\mathrm{Lip}}

\newcommand{\Nbar}{\overline{N}}
\newcommand{\pf}{\par{\noindent\textbf{Proof.~}}}

\begin{document}


\title[Inverse spectral positivity for surfaces]{Inverse spectral positivity for surfaces}

\author[Pierre B\'{e}rard and Philippe Castillon]{Pierre B\'{e}rard and Philippe Castillon}

\date{November 2014}

\begin{abstract}
Let $(M,g)$ be a complete non-compact Riemannian surface. We
consider operators of the form $\Delta + aK + W$, where $\Delta$ is
the non-negative Laplacian, $K$ the Gaussian curvature, $W$ a
locally integrable function, and $a$ a positive real number.
Assuming that the positive part of $W$ is integrable, we address the
question ``What conclusions on $(M,g)$ and on $W$ can one draw from the
fact that the operator $\Delta + aK + W$ is non-negative ?''

\noi As a consequence of our main result, we get a new proof of
Huber's theorem and Cohn-Vossen's inequality, and we improve earlier
results in the particular cases in which $W$ is non-positive and $a
= \4$ or $a \in (0,\4)$.
\end{abstract}

\maketitle

\bigskip

\bigskip

\emph{Published in Rev. Mat. Iberoamericana 30 (2014), 1237-1264.}

\thispagestyle{empty}

\vspace{1.5cm}


\section{Introduction}\label{S-intro}

Let $(M,g)$ be a complete non-compact Riemannian surface. In the
sequel, we will always implicitly assume that $M$ is connected and
orientable, without boundary. We denote by $\Delta$ the non-negative
Laplacian, by $K$ the Gaussian curvature, and by $\mu$ the
Riemannian measure associated with the metric $g$. We denote by $\chi(M)$
the Euler-Poincar\'{e} characteristic of $M$, with the convention that
$\chi(M) = - \infty$ if $M$ does not have finite topology, see for example \cite[Section~1]{Cas06}. \medskip

In this paper, we consider operators of the form $\Delta + aK + W$,
where $a$ is a positive parameter, and $W$ a locally integrable
function. Such operators appear naturally when one studies minimal
(or constant mean curvature) immersions. Let us mention two
examples. The Jacobi (stability) operator of an isometric minimal
immersion $M \looparrowright \R^3$ into Euclidean $3$-space is
$\Delta + 2K$. More generally \cite[Section~3]{FCSc80}, the Jacobi
operator of a minimal immersion $M \looparrowright \Mh^3$ into a
$3$-manifold with scalar curvature $\Sh$ can be written as $\Delta +
K - (\Sh + \2 |A|^2)$, where $|A|$ is the norm of the second
fundamental form of the immersion.
\medskip

More precisely, this paper is concerned with the following
question:~ \emph{What conclusions on the Riemannian surface $(M,g)$,
and on the function $W$, can one draw from the fact that the
operator $\Delta + a K +W$ is non-negative on $(M,g)$~?} ~\ie from
the fact that the associated quadratic form is non-negative on
Lipschitz functions with compact support in $M$ (or equivalently on
$C^1$-functions with compact support),
\begin{equation}\label{E-star}
0 \le \int_M \big( |df|^2 + a K f^2 + W f^2 \big) \, d\mu\,,
\hspace{1cm} \forall f \in \lip_0(M).
\end{equation}

Before stating our results, we recall some definitions.\medskip

Given a function $W$, we let $W_{+}$ and $W_{-}$ denote respectively
the positive and negative parts of $W$, $W_+ = \max\{W,0\}$ and $W_-
= \max\{-W,0\}$, so that $W = W_+ - W_-$ and $|W| = W_+ +
W_-$.\medskip

\textbf{Definitions}. Let $x \in M$, and let $V(r)$ denote the
volume of the geodesic ball $B(x,r)$ for the metric $g$. We say that
$(M,g)$ has \emph{subexponential volume growth} if
$$
\limsup_{r \to \infty} \frac{\ln V(r)}{r} = 0.
$$
We say that $(M,g)$ has \emph{polynomial volume growth of degree at
most $k$} if
$$
\limsup_{r \to \infty} \frac{V(r)}{r^k} < \infty.
$$
We say that $(M,g)$ has $k$-\emph{subpolynomial volume growth} if
$$
\limsup_{r \to \infty} \frac{V(r)}{r^k} = 0.
$$
These definitions do not depend on the choice of the point $x \in
M$. Note that if $(M,g)$ has polynomial volume growth of degree at
most $k$, then it has $k'$-subpolynomial volume growth for any $k' >
k$.\medskip

Our main result is the following.

\begin{thm}\label{T-GG}
Let $(M,g)$ be a complete non-compact Riemannian surface, and let
$W$ be a locally integrable function on $M$, with $W_+$ integrable.
Assume that the operator $\Delta + a K + W$ is non-negative on $M,$
and that either
\begin{itemize}
    \item[(i)] $a \in (\4, \infty)$, or
    \item[(ii)] $a=\4$, and $(M,g)$ has subexponential volume growth, or
    \item[(iii)] $a \in (0,\4)$, and $(M,g)$ has $k_a$-subpolynomial volume
    growth, with $k_a = 2 + \frac{4a}{1-4a}.$
\end{itemize}
Then,
\begin{itemize}
    \item[(A)] The surface $(M,g)$ has finite topology, and at most quadratic
    volume growth. In particular, $(M,g)$ is conformally equivalent to
    a closed Riemannian surface with finitely many points removed.
    \item[(B)] The function $W$ is integrable on $(M,g)$, and
    $$
    0 \le 2\pi a \, \chi(M) + \int_M W \, d\mu .
    $$
    \item[(C)] If $2\pi a \, \chi(M) + \int_M W \, d\mu = 0$, then $(M,g)$ has
    subquadratic volume growth, and $a K + W \equiv 0$ a.e. on the surface $M$.
\end{itemize}
\end{thm}

\newpage
\textbf{Remarks.}
\begin{enumerate}
    \item Theorem~\ref{T-GG} does not mention the case in which $M$
    is closed. When this is the case, \eqref{E-star} implies that,
    $$
    0 \le 2\pi a \chi(M) + \int_M W\, d\mu,
    $$
    with equality if and only if $aK+W \equiv 0$ on $M$. Indeed \cite{FCSc80},
    it suffices to plug the constant function $\mathbf{1}$ into
    \eqref{E-star}, and to notice that equality holds if and only if $\mathbf{1}$
    is in the kernel of the operator $\Delta + aK + W$. When $W_{+} \equiv 0$
    and $W_{-} \not \equiv 0$, the preceding inequality implies that $(M,g)$ is
    conformally equivalent to the round sphere; when $W_{+} \equiv 0$ and $\chi(M)=0$,
    it implies that $W \equiv 0$ and $K \equiv 0$, \ie that $(M,g)$ is a flat torus.
    \item The volume growth assumptions (ii) and (iii) in the statement of the theorem
    are sharp, see Section~\ref{SS-ptgg3} for more details.
\end{enumerate}\medskip

As a corollary of Assertions~A and B in Theorem~\ref{T-GG}, we
obtain Huber's theorem and Cohn-Vossen's inequality.

\begin{thm*}[Huber's theorem]\label{T-Huber}
Let $(M,g)$ be a complete Riemannian surface. Assume that the
negative part $K_{-}$ of the Gaussian curvature is integrable on
$M$. Then,
\begin{itemize}
    \item[(i)] $(M,g)$ has finite topology and is conformally equivalent to a
    closed Riemannian surface with at most finitely many points
    removed;
    \item[(ii)] the Gaussian curvature $K$ is integrable, and $\int_M
    K\, d\mu \le 2\pi \chi(M)$ (Cohn-Vossen's inequality).
\end{itemize}
\end{thm*}

\pf Indeed, write $\Delta = \Delta + K + W$ with $W = - K$, so that
$W_{+} = K_{-}$, and apply Assertions~A and B of Theorem~\ref{T-GG},
Case (i).\qed \medskip

The case in which $W_{+} \equiv 0$ is of particular interest. The
next two results follow directly from Theorem~\ref{T-GG} and its
proof.

\begin{thm}\label{T-G}
Let $(M,g)$ be a complete non-compact Riemannian surface, and let
$q$ be a non-negative locally integrable function on $M.$ Assume
that the operator $\Delta + a K - q$ is non-negative on $M,$ and
that either
\begin{itemize}
    \item[(i)] $a \in (\4, \infty)$, or
    \item[(ii)] $a=\4$, and $(M,g)$ has subexponential volume growth, or
    \item[(iii)] $a \in (0,\4)$, and $(M,g)$ has $k_a$-subpolynomial volume
    growth, with $k_a = 2 + \frac{4a}{1-4a}.$
\end{itemize}
Then,
\begin{itemize}
    \item[(A)] The surface $(M,g)$ has at most quadratic volume growth and
    is conformally equivalent to $\C$ or to $\Cb$ with the standard metrics.
    \item[(B)] The function $q$ is integrable on $(M,g)$, and $\int_M
    q \, d\mu \le 2\pi a \, \chi(M)$.
    \item[(C)] If $M$ is a cylinder, then $(M,g)$ has at most linear volume
    growth and $q \equiv 0$.
\end{itemize}
\end{thm}

\begin{thm}\label{T-C}
Let $(M,g)$ be a complete non-compact $2$-dimensional cylinder.
Assume that the operator $\Delta + a K$ is non-negative on $(M,g)$,
and that either,
\begin{itemize}
    \item[(i)] $a > \4$, or
    \item[(ii)] $a = \4$, and $(M,g)$ has subexponential volume growth, or
    \item[(iii)] $a \in (0,\4)$, and $(M,g)$ has $k_a$-subpolynomial volume
    growth, with $k_a = 2 + \frac{4a}{1-4a}$.
\end{itemize}
Then, $(M,g)$ is flat, \ie $K \equiv 0$.
\end{thm}

\textbf{Remarks.}~\begin{enumerate}
    \item Theorem~\ref{T-GG} and \ref{T-G} extend to the case in which the
    operator $\Delta + a K + W$ is only assumed to have finite
    index. The conclusions, under the assumptions of the theorems, are
    that $(M,g)$ is conformally equivalent to a closed Riemannian surface with a
    finite number of points removed, and that $W$ is integrable over
    $(M,g)$. We refer to Section~\ref{S-fio} for a precise statement
    and its proof.
    \item Another interesting situation, with applications to minimal
    and \textsc{cmc} surfaces, occurs when the potential $q$
    has a positive lower bound, \ie $c = \inf q > 0$. Theorem~\ref{T-G}
    can be extended to this situation, improving the results of
    \cite{Esp09}. We refer to \cite{BerCas12a} for the details.
\end{enumerate} \medskip

Our next result provides an intrinsic version of the optimal length
estimate of L.~Mazet \cite{Maz09}. Note that this is a local result,
we do not need $M$ to be complete. It applies when $M$ is a stable
constant mean curvature surface, possibly with boundary,
isometrically immersed into a simply connected space form $\Mh$, see
Corollary~\ref{C-R}. Our proof follows the same ideas as in Mazet's
paper. We clarify the argument by applying a transplantation method.

\begin{thm}\label{T-R}
Let $(M,g)$ be a Riemannian surface (possibly with boundary
$\partial M$). Assume that the Gaussian curvature satisfies $K \le
\alpha^2$ for some $\alpha > 0$. Let $J$ be the operator $J = \Delta
+ a K - c$, with $a \in [\2, \infty)$ and $c \ge (a + 2) \alpha^2$.
\begin{itemize}
\item[(i)] If $J$  is non-negative in a geodesic ball $B(x,R)$ contained in
$M\setminus\partial M$, then $R \le \frac{\pi}{2\alpha}$.
\item[(ii)] Assume that the geodesic ball $B(x, \frac{\pi}{2\alpha})$ is contained in
$M\setminus\partial M$. If $J$ is non-negative in this ball, then
$c=(a+2)\alpha^2$, $K \equiv \alpha^2$, and $B(x,
\frac{\pi}{2\alpha})$ is covered by the hemisphere
$S_{+}^2(\alpha^2)$ in the sphere with constant curvature
$\alpha^2$.
\end{itemize}
\end{thm}

\textbf{Remark.} In the proof of Theorem~\ref{T-R}, we use the
following classical result. Let $\rho : (\Mh, \gh) \to (M,g)$ be a
Riemannian covering. Let $V$ be a locally integrable function on
$M$, and let $\hat{V} = V \circ \rho$. According to \cite[Theorem~1]{FCSc80},
$\Delta + V \ge 0$ on $(M,g)$ implies that $\hat{\Delta}
+ \hat{V} \ge 0$ on $(\Mh, \gh)$. It is a natural question to
investigate under which conditions the converse statement holds. A
partial answer is given by \cite[Proposition~2.5]{MePeRo08}. As a
matter of fact, one can show that the converse holds
provided that the group $\pi_1(\Mh)$ is co-amenable in the group
$\pi_1(M)$. We defer the precise statement and its proof to
\cite{BerCas12b} because they rely on techniques and ideas different
from those used in the present paper. \bigskip

\textbf{Some background for Theorems~\ref{T-GG}, \ref{T-G} and
\ref{T-C}.}\smallskip

\noindent \textbf{1.} The idea behind the proof of Theorem~\ref{T-GG} goes
back to A.~Pogorelov's proof \cite{Pog81} that orientable
stable minimal surfaces in $\R^3$ are planes. For this purpose, he
shows that a complete simply-connected surface, with non-positive
curvature and non-negative operator $\Delta + 2K$, must be
pa\-ra\-bo\-lic. Another proof consists in showing that such a
simply-connected surface cannot be conformally equivalent to the
unit disk, \cite{CaPe79, FCSc80}. In the latter paper, D.~Fischer-Colbrie
and R.~Schoen prove that there exists
no complete metric $g$ on the unit disk $\D$, such that $\Delta_g + a K_g
\ge 0$ for some $a \ge 1$. More precisely, they show that the set
$I(\D,g) = \{a \ge 0 ~|~ \Delta_g + a K_g \ge 0 \}$ is a closed
interval which does not contain $1$, and they ask what is the value
of the supremum of $I(\D,g)$. This question motivated \cite{Cas06}
and the present paper.\smallskip

\noindent \textbf{2.} A. Pogorelov's result was extended to the case in which
$\Delta + a K \ge 0$ for some $a > \4$ by S.~Kawai \cite{Kaw88}. A
more general setting (general topology and curvature, $a > \2$) was
considered by R. Gulliver and B. Lawson in \cite{GuLa86}.
A.~Pogorelov's method was improved by T.~Colding and W.~Minicozzi
\cite{CoMi02} and, later on, by Ph.~Castillon \cite{Cas06} who first
proved Case (i) in Theorem~\ref{T-G} (with $q \equiv 0$); see also
\cite{EsRo10} and \cite{MePeRo08} which contain applications to
constant mean curvature surfaces in $3$-manifolds. Cases (ii) and
(iii) in Theorem~\ref{T-G} were first considered by J.~Espinar and
H.~Rosenberg in \cite{EsRo10}, under more restrictive assumptions on
$(M,g)$. \smallskip

\noindent \textbf{3.} The main new idea in the proof of Theorem~\ref{T-GG} is
to introduce the function $\chih(t) := \sup \{ \chi\big( B(s) \big)
~|~ s \in [t, \infty) \}$, the supremum of the Euler-Poincar\'{e}
characteristics of open geodesic balls with radius at least $t$,
whose jumps describe the large scale topology of $M$, see
Section~\ref{S-preli}. We also introduce new functions to test the
positivity of the quadratic form \eqref{E-star}, see
Lemmas~\ref{L-A} and \ref{L-B}.\smallskip

\noindent \textbf{4.} Theorem~\ref{T-GG}, Assertion C, and Theorem~\ref{T-C}
were motivated by the following result due to D.~Fi\-scher-Colbrie
and R.~Schoen \cite[Theorem~3]{FCSc80} reformulated for our
purpose.

\begin{thm*}
Let $N$ be a complete oriented $3$-manifold of non-negative scalar
curvature. Let $M$ be an oriented complete non-compact stable
minimal surface in $N$. Then $M$ is conformally equivalent to the
complex plane or to a cylinder. If $M$ is a cylinder, and the
absolute total curvature of $M$ is finite, then $M$ is flat and
totally geodesic.
\end{thm*}

The proof involves the operator $\Delta + K$. D.~Fischer-Colbrie and
R.~Schoen also point out (Remark 2, p. 207) that ``the assumption of
finite total curvature should not be essential''. This is indeed the
case: see \cite[proof of Theorem~2]{ScYa82}, and \cite{Miy93} for a
proof using $L^2$ harmonic $1$-forms. Theorem~\ref{T-G} and
Theorem~\ref{T-C}, with $a=1$, provide another method to ans\-wer
the question raised by D.~Fischer-Colbrie and R.~Schoen
positively.\medskip

\textbf{Remark.} Case~(i) in Theorem~\ref{T-C} appears in
\cite[Section~3.3]{Rei10}  under the assumption that $a \ge 1$; in
\cite[Theorem~6.3]{EsRo10} under the restrictive assumption that
$\int_M K^{+} \, d\mu$ is finite; in \cite{Esp11} with a different
proof. \bigskip

The paper is organized as follows. In Section~\ref{S-preli}, we fix
the notations and state some technical lemmas to be used later on.
The proof of Theorem~\ref{T-GG} is given in Section~\ref{SS-ptgg1}
(Assertions~A and B) and Section~\ref{SS-ptgg2} (Assertion~C). The
fact that the volume growth assumptions in the theorem are optimal
is explained in Section~\ref{SS-ptgg3}. In Section~\ref{S-fio} we
give extensions of Theorems~\ref{T-GG} and \ref{T-G} to the case of
finite index operators. The proofs of Theorems~\ref{T-G}, \ref{T-C}
and \ref{T-R} are given in the last section.\medskip

The first author was partially supported by \textsc{cnrs} and
by the cooperation programmes \textsc{arcus} Rh\^{o}ne-Alpes--Br\'{e}sil and
\textsc{math-AmSud} during the preparation of this paper.
The authors would like to thank J.~Espinar and L.~Mazet for their
comments on a preliminary version.


\section{Notations and preliminary results}\label{S-preli}

In this section, we fix some notations which will be used throughout
the paper, and we state some preliminary results.

\subsection{Notations}\label{SS-not} In this paper, unless otherwise stated,
$(M,g)$ denotes a complete non-compact Riemannian surface without
boundary. We also assume that $M$ is connected and orientable.\smallskip

\noindent \textbf{1.} The non-negative Laplacian for the metric $g$ will be
denoted by $\Delta$, the Gauss curvature by $K$ and the Riemannian
measure by $\mu$. \smallskip

\noindent \textbf{2.} Let $x_0$ be a given reference point
in $M$. We let $r(x)$ denote the Riemannian distance from the point
$x$ to the point $x_0$.  We let $B(s)$ denote the \emph{open}
geodesic ball with center $x_0$ and radius $s$. For $t < s$, we let
$C(t,s)$ denote the open set $C(t,s) =
B(s)\!\setminus\!\overline{B}(t)$. The volume of the ball $B(s)$ is
denoted by $V(s)$, the length of the boundary of $B(s)$ by $L(s)$.
The length function is a priori only defined for $s \in
\R_+\!\setminus\!E$, where the set of exceptional values $E$ is
closed, of Lebesgue measure zero. On the set $\R_+\!\setminus\!E$,
the function $L$ is $C^1$ and satisfies Fiala's inequality, see
\cite{Fia40},
\begin{equation}\label{E-fiala}
L'(t) \le 2\pi \chi \big( B(t) \big) - \int_{B(t)} K \, d\mu,
\end{equation}
where $\chi \big( B(t) \big)$ is the Euler-Poincar\'{e} characteristic
of the ball $B(t)$. As a matter of fact, the function $L$ can be
extended to $\R_{+}$. This follows from the work of F.~Fiala
\cite{Fia40}, P.~Hartman \cite{Har64}, and K.~Shiohama and M.~Tanaka
\cite{ShTa89, ShTa93}. More precisely, there exist two real
functions $H, J$ defined on $\R_+$, with $H$ absolutely continuous
on any compact subset, and $J$ non-decreasing, such that $H(s)-J(s)$
coincides with $L(s)$ when $s$ is not in $E$. The set $E$ and the
function $J$ are defined in terms of the cut locus of the point
$x_0$. The (extended) function $L$ is not continuous in general
\cite[Figure~1]{Har64}. However, it satisfies,
\begin{equation}\label{E-lf0}
L(t^+) = L(t) \text{~and~} L(t^-) \ge L(t), ~~\forall t > 0.
\end{equation}
Furthermore, the function $V$ is differentiable almost everywhere,
and $V'(s) = L(s)$. From the formula $L=H-J$, one can deduce that
\begin{equation}\label{E-lf1}
L(b) - L(a) \le L(b^-) - L(a) \le \int_a^b L'(t) \, dt, \text{
~whenever~ } 0 \le a < b.
\end{equation}

\textbf{Remark}. In Fiala's paper, $M=\R^2$ and $g$ is real
analytic. In this case, the set $E$ is discrete. Hartman's paper
considers the case $(\R^2,g)$, with $g$ smooth enough. The papers of
Shiohama and Tanaka deal with the general case in which $M$ may have
finite or infinite topology. All these papers rely on a sharp
analysis of the cut locus of a simple closed curve, and on the
differential inequality (\ref{E-fiala}) satisfied by the length
function $L$ away from the exceptional set $E$. This was initiated
by Fiala, refined by Hartman and later by Shiohama-Tanaka to take
into account the transitions from a real analytic to a smooth
metric, and from $\R^2$ to a general surface $M$.\smallskip

\noindent \textbf{3.} We introduce the total curvature of the ball
$B(s)$ to be $$G(s) = \int_{B(s)} K(x) \, d\mu(x).$$

\noindent \textbf{4.} We denote by $\chi\big(B(s)\big)$ the
Euler-Poincar\'{e} characteristic of the open ball $B(s)$. We introduce
the function
\begin{equation}\label{E-sep1}
\chih (s) = \sup\{\chi \big( B(t) \big) ~|~ t \in [s, \infty)\}.
\end{equation}
Both functions are continuous to the left. The function $\chih$ is a
non-increasing function from $[0, \infty)$ to $\Z$. It has at most
countably many discontinuities. We write them as a sequence, finite
possibly empty, or infinite tending to infinity,
\begin{equation}\label{E-sep2}
\{t_j\}_{j=1}^{\Nbar} = \{ 0 < t_1 < t_2 < \cdots < t_n < \cdots \},
\end{equation}
with $\Nbar \in \N \cup \{\infty\}$, $\Nbar = 0$ when the sequence
is empty, and $\Nbar = \infty$ when it is infinite. Note that this
sequence depends on the reference point $x_0$. \smallskip

At the discontinuity $t_n, n\ge 1$, the function $\chih$ has a jump
\begin{equation}\label{E-sep3}
\omega_n = \chih (t_n^-) - \chih (t_n^+), \text{~with~} \omega_n \in
\N, \omega_n \ge 1.
\end{equation}
Therefore,
\begin{equation}\label{E-sep4}
\left\{%
\begin{aligned}
    \chih (s) &=1, \text{~for~} s \in [0, t_1], \text{~and}\\
    \chih (s) &= 1 - \big(\omega_1 + \cdots + \omega_n\big) \le - (n-1),
    \text{~for~} s \in (t_n, t_{n+1}].
\end{aligned}%
\right.
\end{equation}\smallskip

The function $\chih$ somehow describes the large scale topology of
$M$ as the following lemma shows.

\begin{lem}\label{L-top}
Let $(M,g)$ be a complete Riemannian surface. Let
$\{t_j\}_{j=1}^{\Nbar}$ be the discontinuities of the function
$\chih$, with jumps $\{\omega_j\}$, relative to some reference point
$x_0$ in $M$. Let $\chi(M)$ be the Euler-Poincar\'{e} characteristic of
$M$, with $\chi(M) = - \infty$ if $M$ does not have finite topology.
Then
$$
1 - \sum_{n=1}^{\Nbar} \omega_n \le \chi(M).
$$
\end{lem}

\pf Apply \cite[Lemma 1.4]{Cas06}.

\noid If $M$ has finite topology, then there exists a value $s_0$
such that $\chi\big(B(s)\big) \le \chi(M)$ for all $s\ge s_0$. By
(\ref{E-sep4}), this implies that $1 - \sum_{n=1}^{\Nbar} \omega_n
\le \chih(s_0) \le \chi(M)$.

\noid Otherwise, $\chi\big( B(s) \big)$ tends to minus infinity when
$s$ tends to infinity, so does $\chih(s)$, and formula
(\ref{E-sep4}) implies that $1 - \sum_{n=1}^{\Nbar} \omega_n = -
\infty$. \qed \smallskip

\noindent \textbf{5.} As mentioned earlier, the Euler-Poincar\'{e} characteristic
of balls is related to the length function and to the total
curvature of balls. More precisely, we have the inequalities,
\begin{equation}\label{E-lec}
\left\{%
\begin{aligned}
    & \text{For all~} 0 \le a < b,\\
    & L(b^{-}) - L(a) \le 2\pi (b-a) \sup\{\chi\big(B(s)\big) ~|~ s \in [a,b]\} -
    \int_a^b G(s) \, ds,\\
    & L(b^{-}) - L(a) \le 2\pi (b-a) \chih(a) - \int_a^b G(s) \, ds,\\
\end{aligned}%
\right.
\end{equation}
which follow by integrating the inequality (\ref{E-fiala}) satisfied
by $L'(t)$ for $t \in ]0,\infty[ \setminus E$ \cite[p.~326-328]{Fia40},
\cite[Proposition~6.1]{Har64}, \cite[Proposition~3.7]{ShTa89}.
Note that we can substitute $L(b^{-})$ by $L(b)$
in (\ref{E-lec}), because of inequality (\ref{E-lf0}).

\subsection{Technical lemmas}\label{SS-tls}

\textbf{Definition}. Let $0 \le R < S$. We say that a function $\xi
: [R,S] \to \R$ is \emph{admissible} in the interval $[R,S]$ if
\begin{equation*}
\left\{%
\begin{aligned}
& \xi \text{~is~} C^1 \text{~and piecewise~} C^2 \text{~in~} [R,S],\\
& \xi \ge 0, \; \xi' \le 0 \text{~and~} \xi'' \ge 0.
\end{aligned}%
\right.
\end{equation*}

The next two lemmas extend Lemma~1.8 in \cite{Cas06}, whose proof
uses the method of \cite{CoMi02}.

\begin{lem}\label{L-t0}
For all $0 \le a < b$, and for all admissible functions $\xi$ on
$[a,b]$,
\begin{equation}\label{E-tl}
\begin{split}
\int_{C(a,b)} K(x) \xi^2\big(r(x)\big) \, d\mu(x) \le & ~\xi^2G
\Big|_a^b
- 2\pi \chih(a) \xi^2 \Big|_a^b + (\xi^2)'L\Big|_a^{b^{-}} \\
&  - \int_{C(a,b)} (\xi^2)'' \big(r(x)\big) \, d\mu(x).
\end{split}%
\end{equation}
\end{lem}

Note that in the right-hand side of equation (\ref{E-tl}) one can substitute
$(\xi^2)'L\big|_a^{b^{-}}$ by $(\xi^2)'L\big|_a^b$, using
(\ref{E-lf0}) and the fact that $\xi'$ is non-positive. \medskip

\pf We sketch the proof for completeness. First assume that $\xi$ is
$C^2$. By the co-area formula,
$$
\int_{C(a,b)} K \xi^2(r) \, d\mu = \int_a^b \xi^2(t) G'(t) \, dt,
$$
where $G(t)$ is the total curvature of the ball $B(t)$. Introduce
the function $H(t) := \int_a ^t G(s) \, ds$, and integrate the
preceding equality by parts twice to get,
$$
\int_{C(a,b)} K \xi^2(r) \, d\mu = \xi^2 G\Big|_a^b -
(\xi^2)'(b)H(b) + \int_a^b H(t) (\xi^2)''(t) \, dt.
$$
One can estimate $H$ in the right-hand side using (\ref{E-lec}) and
the signs of $\xi$ and its derivatives. After some computations and
applying the co-area formula once more, one obtains,
$$
\int_{C(a,b)} K \xi^2(r) \, d\mu \le \big\{ \xi^2G - 2\pi \chih(a)
\xi^2 + (\xi^2)'L \big\} \Big|_a^{b^-} - \int_{C(a,b)} (\xi^2)''(r)
\, d\mu.
$$

This proves the lemma when $\xi$ is $C^2$. The fact that the lemma
holds for $C^1$ and piecewise $C^2$ functions $\xi$ follows by
cutting the interval into subintervals in which $\xi$ is $C^2$.
Apply the preceding method in each sub-interval $(c,d) \subset
(a,b)$ using an inequality similar to (\ref{E-lec}) with $\chih(a)$
in place of $\chih(c)$; use the fact that $\xi$ and $\chih$ are
non-increasing, and the inequality $L(t^-) \ge L(t^+) = L(t)$ to
conclude. \qed \medskip

Taking into account the discontinuities  $\{t_n\}_{n\ge 1}$ of the
function $\chih$, see Section~\ref{SS-not}, formula~(\ref{E-sep4}),
we have the following lemma.

\begin{lem}\label{L-t1}
Let $\{t_j\}_{j=1}^{\Nbar}$ be the discontinuities of the function
$\chih$. Define the index $N(R)$ to be the largest integer $n$ such
that $t_{n} \le R$. Let $t_0 = 0$.

Let $\xi$ be an admissible function in the interval $[R,Q]$. Then,
\begin{equation}\label{E-T3}
\begin{split}
& \int_{C(R,Q)} K \xi^2(r) \, d\mu \le 2\pi
\Big[ \xi^2(R) \chih (t_{N(R)}) - \xi^2(Q) \chih (t_{N(Q)}) \\
& - \sum_{n=N(R)+1}^{N(Q)} \omega_n \xi^2(t_n) \Big]
+ \xi^2G \Big|_{R}^{Q} + (\xi^2)' L \Big|_{R}^{Q} - \int_{C(R,Q)}
(\xi^2)''(r) \, d\mu.\\
\end{split}
\end{equation}
Taking $R=0$ and assuming that $\xi(Q)=0$, we have the inequality
\begin{equation}\label{E-T0}
\int_{B(Q)} K \xi^2(r) \, d\mu \le 2\pi
\left\{\xi^2(0) - \sum_{n=1}^{N(Q)} \omega_n \xi^2(t_n) \right\}
- \int_{B(Q)} (\xi^2)''(r) \, d\mu.
\end{equation}

In particular, assuming that $\xi(Q)=0$, we have the inequality
\begin{equation}\label{E-T1}
\int_{B(Q)} K(x) \xi^2(r) \, d\mu \le 2\pi \xi^2(0) - \int_{B(Q)}
(\xi^2)''(r) \, d\mu.
\end{equation}
\end{lem}

\pf To prove (\ref{E-T3}), split the integral $\int_{C(R,Q)} K
\xi^2(r) \, d\mu$ into a sum,
$$\int_{C(R,Q)} = \int_{C(R,t_{N(R)+1})} +
\sum_{n=N(R)+1}^{N(Q)-1} \int_{C(t_n,t_{n+1})} +
\int_{C(t_{N(Q)},Q)},$$
apply Lemma~\ref{L-t0} and use (\ref{E-sep4}). To establish the last
two inequalities, use the fact that $\xi(Q)=0$ and $G(0)=L(0)=0$.
\qed\medskip

The next two lemmas provide admissible functions which we will plug
into \eqref{E-star} later on.

\begin{lem}\label{L-A}
Fix $0 < R < 5R < Q$, and define the function $\xi_{\alpha, \beta,
R, Q}$ by
\begin{equation*}\label{E-A1}
\xi_{\alpha, \beta, R, Q}(t) = \left\{%
\begin{aligned}
& e^{(1-\frac{t}{2R})^2}, \text{~for~} 0 \le t \le R,\\
& \beta \big( e^{-\alpha t} - e^{-\alpha Q}\big), \text{~for~} R \le t
\le Q.\\
\end{aligned}
\right.
\end{equation*}
Then, there exists a unique choice $\alpha(R,Q)$, $\beta(R,Q)$ of
the parameters $\alpha, \beta$ such that the corresponding function
$\xi_{R,Q}$ is admissible in the interval $[0,Q]$. Furthermore,
$$
1 \le 4 R \, \alpha(R,Q) \le 2 \text{~and~} 1 \le \beta(R,Q) \le 10.
$$
\end{lem}

\begin{lem}\label{L-B}
For $a \in (0, \4)$, let $\alpha = \frac{2a}{1-4a}$ and $\beta =
\frac{a}{1-4a}$. For $0 < R < Q$ and $0 < \delta, \epsilon$, let
$\xi_{\delta, \epsilon,R,Q}$ be the function,
\begin{equation}\label{E-B4}
\xi(t) = \left\{%
\begin{aligned}
    & (1 + \frac{t}{R})^{-\beta}, \text{~for~} t \in [0,R],\\
    & \delta \Big( (1 + \epsilon t)^{-\alpha} - (1 + \epsilon Q)^{-\alpha} \Big),
    \text{~for~} t \in [R,Q].\\
\end{aligned}%
\right.
\end{equation}
There exists a positive constant $C(\alpha,\beta) > 1$, such that
for $0 < R \le C(\alpha,\beta)R < Q$, there is a unique choice
$\delta(R,Q)$ and $\epsilon(R,Q)$ of the parameters $\delta,
\epsilon$, such that the function $\xi_{R,Q}$ defined by equation
(\ref{E-B4}) is admissible in the interval $[0,Q]$. Furthermore,
there exist positive constants $c_1, c_2$ such that
$$1 \le 6 R \, \epsilon (R,Q) \le 2 \text{ ~and~ }
c_1 \le \delta(R,Q) \le c_2.$$
\end{lem}\smallskip

We leave the proofs of these lemmas to the reader.

\section{Proof and optimality of Theorem~\ref{T-GG}}\label{S-ptgg}

\subsection{Proof of Theorem~\ref{T-GG} -- Assertions A and B}\label{SS-ptgg1}

In Step 1, we make some preparation. In Step~2, we prove that $M$
has finite topology, and that $W$ is integrable and satisfies $0 \le
2 \pi a \chi(M) + \int_M W \, d\mu$. We can actually finish the
proof of assertions A and B in the Case (i). In Step~3, we prove
that $(M,g)$ has at most quadratic volume growth. Steps 2 and 3 both
follow from adequate choices of test functions (using
Lemma~\ref{L-A} and \ref{L-B}) in the stability condition
\eqref{E-star}, depending on the case at hand (i), (ii) or (iii).

\subsubsection{Step 1} We choose an admissible function $\xi$ on
$[0,Q]$, with $\xi(Q)=0$, and we apply the stability condition
\eqref{E-star} to the Lipschitz function $\xi(r)$, where $r$ is the
Riemannian distance to some given point $x_0 \in M$. We obtain,
$$
\int_{B(Q)} W_-\, \xi^2(r) \, d\mu \le \int_{B(Q)} \big\{
(\xi')^2(r) + a K \xi^2(r) \big\} \, d\mu + \int_{B(Q)} W_+\,
\xi^2(r) \, d\mu. 
$$

Because $\xi$ is admissible in $[0,Q]$ and $\xi(Q)=0$, we can apply
Lemma~\ref{L-t1}, inequality (\ref{E-T0}), and we obtain
\begin{equation}\label{E-b}
\begin{split}
& \int_{B(Q)} W_{-}\, \xi^2(r) \, d\mu +  2 \pi a  \sum_{n=1}^{N(Q)}
\omega_n \xi^2(t_n) \le 2 \pi a \xi^2(0) \\
& + \int_{B(Q)} \big\{ (1-2a) (\xi')^2(r) - 2a (\xi
\xi'')(r) \big\} \, d\mu
+ \int_{B(Q)} W_{+}\, \xi^2(r) \, d\mu,
\end{split}
\end{equation}
where we have used the notations of Lemma~\ref{L-t1}. Inequality
\eqref{E-b} holds for all admissible functions $\xi$ in $[0,Q]$
which vanish at $Q$. \medskip

Recall (Section~\ref{SS-not}, \S 4) that the points of discontinuity of
the function $\chih$ form a sequence $\{t_n\}_{n=1}^{\Nbar}$ which
is either finite possibly empty, or infinite tending to infinity,
with stopping index $\Nbar \in \N \cup \{\infty\}$. \medskip

We fix $N$ to be either the stopping index $\Nbar$, if $\Nbar \in
\N$, or any fixed integer otherwise. We also fix some $R$, with $0 <
R < Q$. For $Q$ large enough, $Q > t_N$ and $Q \ge C(\xi) R$,
inequality \eqref{E-b} implies that
\begin{equation}\label{E-c}
\begin{split}
& \int_{B(R)} W_{-}\, \xi^2(r) \, d\mu + 2 \pi a  \sum_{n=1}^{N}
\omega_n \xi^2(t_n) \le 2 \pi a \xi^2(0) \\
& + \int_{B(Q)} \big\{ (1-2a) (\xi')^2(r) - 2a (\xi
\xi'')(r) \big\} \, d\mu
+ \int_{B(Q)} W_{+}\, \xi^2(r) \, d\mu,
\end{split}
\end{equation}
where this inequality holds for any admissible function $\xi$ in
$[0,Q]$ vanishing at $Q$, and for any fixed $N$ and $R$ as above.
\medskip

The idea is now to apply \eqref{E-c} to a function $\xi$ which is
well adapted to the case at hand, (i), (ii) or (iii), and to the
assertion we want to prove.

\subsubsection{Step 2} We will now show that $M$ has finite
topology, and that $W$ in integrable over $(M,g)$. We consider the
cases (i), (ii) and (iii) separately. \medskip

\boxit{\hbox{Case (i).}} Here, $a \in (\4, \infty)$. Choose $\xi (t)
= (1 - \frac{t}{Q})^{\alpha}$ for $t \in [0,Q]$, with $\alpha \ge
1$. Then,
$$
(1-2a) (\xi')^2 - 2a \xi \xi'' = - \frac{\alpha [(4a-1)\alpha -
2a]}{Q^2} (1 - \frac{t}{Q})^{2\alpha-2}. 
$$
We now fix some $\alpha > \frac{2a}{4a-1}$. Plugging the previous
equality into inequality \eqref{E-c} yields that, for all $R$ and
$N$ fixed,
\begin{equation}\label{E-e1}
\begin{split}
& \int_{B(Q)} W_{-}\, \, \xi^2(r) \, d\mu + 2 \pi a \sum_{n=1}^{N}
\omega_n \xi^2(t_n) \\
& + \frac{\alpha [(4a-1)\alpha - 2a]}{Q^2} \int_{B(Q)}
(1 - \frac{r}{Q})^{2\alpha-2} \, d\mu \le 2 \pi a + \int_{B(Q)}
W_{+}\, \xi^2(r) \, d\mu.
\end{split}
\end{equation}

Note that the three terms in the left-hand side of \eqref{E-e1} are
non-negative. Using the fact that $W_{+}$ is integrable, we obtain
that
$$
\sum_{n=1}^N \omega_n (1 - \frac{t_n}{Q})^{2\alpha} \le 1 +
\frac{1}{2\pi a} \int_{M} W_{+} \, d\mu.
$$
Letting $Q$ tend to infinity, we conclude that
$$
\sum_{n=1}^N \omega_n  \le 1 + \frac{1}{2\pi a} \int_{M} W_{+} \,
d\mu,
$$
for any fixed $N$ as above. It follows that $\Nbar$ is actually
finite and, by Lemma~\ref{L-top}, that $M$ has finite topology.

\begin{rem}\label{R-1}
When $W_{+} \equiv 0$, the preceding inequality implies that
$\Nbar=0$, in which case $M$ is homeomorphic to $\C$, or that $\Nbar
= 1$ and $\omega_1=1$, in which case $M$ is homeomorphic to $\C$ or
to $\Cb$.
\end{rem}

From \eqref{E-e1} and the previous conclusions, we can choose $N
=\Nbar$ and we obtain,
$$
(1-\frac{R}{Q})^{2\alpha} \int_{B(R)} W_{-}\, d\mu \le 2\pi a \big(
1 - \sum_{n=1}^{\Nbar} \omega_n \xi^2(t_n)\big) + \int_{M} W_{+} \,
d\mu.
$$
Letting $Q$ tend to infinity and using Lemma~\ref{L-top}, this
proves that
$$\int_{B(R)} W_{-}\, d\mu \le 2\pi a \chi(M) + \int_{M}
W_{+} \, d\mu.
$$
Since this is true for any $R>0$, we have that $W_{-}$ is
integrable, and Assertion~B follows in the Case (i).\medskip

Note that in Case~(i), we can also conclude that $(M,g)$ has at most
quadradic volume growth. Indeed, from \eqref{E-e1}, we can infer
that there exists a positive constant $C_{\alpha}$ such that
$$
C_{\alpha} Q^{-2} V(\frac{Q}{2}) \le 2\pi a \chi(M) + \int_{M} W_{+}
\, d\mu.
$$

\begin{rem}\label{R-2}
When $W_{+} \equiv 0$, we actually get a sharper result when $M$ is
homeomorphic to $\Cb$ (\ie to a cylinder). Indeed, in that case
$\Nbar = 1$ and $\omega_1 = 1$, and inequality \eqref{E-e1} gives
that
$$
C_{\alpha} Q^{-2} V(\frac{Q}{2}) \le 2 \pi a \, \big\{ 1 - (1 -
\frac{t_1}{Q})^{2\alpha} \big\}.
$$
It follows that $M$ has at most linear volume growth in this
particular case.
\end{rem}

Note that this completes the proof of Theorem~\ref{T-GG},
Assertions~A and B, in the Case (i). In the following arguments, we
will concentrate on the cases (ii) and (iii). \medskip

\boxit{\hbox{Case (ii).}} Here $a = \4$ and $(M,g)$ has
subexponential volume growth. We choose $\xi(t) = e^{-\alpha t} -
e^{- \alpha Q}$ in $[0,Q]$ for some $\alpha > 0$. Then,
$$
(1-2a) (\xi')^2 - 2a \xi \xi'' = \2 \alpha^2 e^{-\alpha Q}
e^{-\alpha t}. 
$$
Plugging this equality into \eqref{E-c}, we obtain, for all $R$ and
$N$ fixed,
\begin{equation}\label{E-e2}
\begin{split}
& \int_{B(R)} W_{-}\, \xi^2(r) \, d\mu + \frac{\pi}{2} \sum_{n=1}^{N}
\omega_n \xi^2(t_n) \le \frac{\pi}{2} \xi^2(0) \\
& + \2 \alpha^2 e^{-\alpha Q} \int_{B(Q)} e^{-\alpha r}\,
d\mu + \int_{B(Q)} W_{+}\, \xi^2(r) \, d\mu.
\end{split}%
\end{equation}

We have the following lemma,

\begin{lem}\label{L-N1}
If $(M,g)$ has subexponential volume growth, then for any positive
$\alpha$,
$$
\lim_{Q \to \infty} e^{-\alpha Q} \int_{B(Q)} e^{-\alpha r}\, d\mu =
0.
$$
\end{lem}

\pf Use the co-area formula and integration by parts. \qed
\medskip

Let $Q$ tend to infinity in \eqref{E-e2}, and use Lemma~\ref{L-N1}
to obtain,
$$
\int_{B(R)} W_{-}\, e^{-\alpha r}\, d\mu + \frac{\pi}{2}
\sum_{n=1}^{N} \omega_n e^{-2\alpha t_n} \le \frac{\pi}{2}+ \int_{M}
W_{+} \, d\mu, 
$$
and this inequality holds for all $\alpha > 0$ and $N, R$ as above.
Letting $\alpha$ tend to zero, we can conclude as in Case (i) that
$\Nbar$ is finite and hence that $M$ has finite topology, and that
$W$ is integrable and satisfies,
$$
0 \le \frac{\pi}{2} \chi(M) + \int_{M} W \, d\mu.
$$

\begin{rem}\label{R-3}
When $W_{+} \equiv 0$, we can conclude as in Remark~\ref{R-1} that
$M$ is homeomorphic to $\C$ or to $\Cb$.
\end{rem}

Note that, unlike in Case (i), we have not yet obtained quadratic
volume growth (see Step 3). \bigskip

\boxit{\hbox{Case (iii).}} Here $a \in (0,\4)$ and $(M,g)$ has
$k_a$-subpolynomial volume growth, with $k_a = 2 + \frac{4a}{1-4a}$.
We choose $\xi(t) = (1 + \epsilon t)^{-\alpha} - (1 + \epsilon
Q)^{-\alpha}$ in $[0,Q]$, with $\alpha = \frac{2a}{1-4a}$ and some
$\epsilon > 0$. Then,
$$
(1-2a) (\xi')^2 - 2a \xi \xi'' = 2a \alpha (\alpha + 1) \epsilon^2
(1 + \epsilon  Q)^{-\alpha}(1 + \epsilon t)^{-\alpha -2}.
$$
Plugging this equality into \eqref{E-c}, we obtain,
\begin{equation}\label{E-e3}
\begin{split}
& \int_{B(R)} W_{-}\, \xi^2(r) \, d\mu  + 2 \pi a \sum_{n=1}^{N}
\omega_n \xi^2(t_n) \le 2 \pi a \xi^2(0) \\
& + 2a \alpha (\alpha + 1) \epsilon^2 (1 + \epsilon
Q)^{-\alpha} \int_{B(Q)} (1 + \epsilon r)^{-\alpha -2} \, d\mu +
\int_{M} W_{+} \, d\mu.
\end{split}
\end{equation}

We have the following lemma,

\begin{lem}\label{L-N2}
Let $(M,g)$ be a Riemannian surface with $k_a$-subpolynomial volume
growth, with $k_a = 2 + \frac{4a}{1-4a}$. Then, for $\alpha =
\frac{2a}{1-4a}$ and any $\epsilon
> 0$,
$$
\lim_{Q \to \infty} (1 + \epsilon  Q)^{-\alpha} \int_{B(Q)} (1 +
\epsilon r)^{- \alpha - 2} \, d\mu = 0.
$$
\end{lem}

\pf Use the co-area formula and integration by parts. \qed.\medskip

Since both terms in the left-hand side of \eqref{E-e3} are
non-negative, letting $Q$ tend to infinity and using
Lemma~\ref{L-N2}, we obtain,
$$
\int_{B(R)} W_{-}\, (1 + \epsilon r)^{-2\alpha} \, d\mu + 2 \pi a
\sum_{n=1}^N \omega_n (1 + \epsilon t_n)^{-2\alpha} \le 2 \pi a +
\int_{M} W_{+} \, d\mu, 
$$
and this inequality holds for any $\epsilon > 0$. Letting $\epsilon$
tend to zero, we obtain,
$$
\int_{B(R)} W_{-}\,  \, d\mu + 2 \pi a \sum_{n=1}^N \omega_n \le 2
\pi a + \int_{M} W_{+} \, d\mu,
$$
and we can conclude as in the previous cases that $\Nbar$ is finite,
that $M$ has finite topology, and that $W$ is integrable, with
$$
0 \le 2\pi a \chi(M) + \int_M W\, d\mu.
$$

\begin{rem}\label{R-4}
When $W_{+} \equiv 0$, we can show as  in Remark~\ref{R-1} that $M$
is homeomorphic to $\C$ or to $\Cb$.
\end{rem}

\subsubsection{Step 3}\label{sss-s3} We now show that $(M,g)$ has at most
quadratic volume growth. We have already dealt with Case (i) in
Step~2. We now consider Cases (ii) and (iii). Recall from Step~2
that $\Nbar$ is finite. \medskip

\boxit{\hbox{Case (ii).}} Here $a=\4$, and $(M,g)$ has
subexponential volume growth. We choose the function $\xi$ to be
$\xi_{R,Q}$ as given by Lemma~\ref{L-A},
\begin{equation*}
\xi (t) = \left\{%
\begin{aligned}
& e^{(1 - \frac{t}{2R})^2}, \text{ ~for~ } t \in [0,R]\\
& \beta (e^{-\alpha t} - e^{-\alpha Q}), \text{ ~for~ } t \in [R,Q],
\end{aligned}%
\right.
\end{equation*}
with $0 < R < 5 Q$ and $\alpha, \beta$ given by the lemma, so that
$\xi$ is admissible in $[0,Q]$ and vanishes at $Q$. We apply
\eqref{E-c} again (making $W_{-} \equiv 0$ which is sufficient for
our estimates). For this purpose, we compute,
\begin{equation*}
(\xi')^2 - \xi \xi'' = \left\{%
\begin{aligned}
& - \frac{1}{2R^2} e^{2(1 - \frac{t}{2R})^2}, \text{ ~for~ } t \in [0,R],\\
& ~\alpha^2 \beta^2 e^{-\alpha Q} e^{-\alpha t} , \text{ ~for~ } t \in
[R,Q],
\end{aligned}%
\right. 
\end{equation*}
and we obtain,
\begin{equation}\label{E-q2}
\begin{split}
& \frac{1}{4R^2} \int_{B(R)} e^{2(1 - \frac{r}{2R})^2} \, d\mu  \le
\frac{\pi}{2}
\big\{ e^2 - \sum_{n=1}^{\Nbar}\omega_n e^{2(1 - \frac{t_n}{2R})^2} \big\}\\
& + \2 \alpha^2 \beta^2 e^{-\alpha Q} \int_{C(R,Q)}
e^{-\alpha r} \, d\mu + \int_{B(Q)} W_{+} \xi^2(r) \, d\mu,
\end{split}
\end{equation}
where we have chosen $R > t_{\Nbar}$. We fix $R > t_{\Nbar}$, and we
let $Q$ tend to infinity, using the facts that $\alpha$ and $\beta$
remain controlled, and that the second term in the right-hand side
of \eqref{E-q2} goes to zero when $Q$ tends to infinity because
$(M,g)$ has subexponential volume growth (Lemmas~\ref{L-A} and
\ref{L-N1}). Finally, we obtain,
$$
R^{-2} V(R) \le C \big\{ 1 - \sum_{n=1}^{\Nbar} \omega_n
e^{-\frac{t_n}{R}(2 - \frac{t_n}{2R})} \big\} + e^2\int_M W_+ \,
d\mu \le C',
$$
for some constant $C'$ independent of $R$. This gives that $M$ has at
most quadratic volume growth.

\begin{rem}\label{R-5}
When $W_{+} \equiv 0$ and $\chi(M) = 0$ (which corresponds to $\Nbar
= 1$ and $\omega_1=1$), the above estimate gives that $M$ has at
most linear volume growth.
\end{rem}

\boxit{\hbox{Case (iii).}} Here, $a \in (0,\4)$ and $(M,g)$ has
$k_a$-subpolynomial volume growth, with $k_a = 2 + \frac{4a}{1-4a}$.
We choose the function $\xi$ to be $\xi_{R,Q}$ as given by
Lemma~\ref{L-B},
\begin{equation*}
\xi (t) = \left\{%
\begin{aligned}
& (1 + \frac{t}{R})^{-\beta}, \text{ ~for~ } t \in [0,R],
\beta = \frac{a}{1-4a},\\
& ~ \delta \big\{ (1 + \epsilon t)^{-\alpha} - (1 + \epsilon
Q)^{-\alpha} \big\}, \text{ ~for~ } t \in [R,Q], \alpha =
\frac{2a}{1-4a},
\end{aligned}%
\right.
\end{equation*}
with $0 < R \ll Q$ and $\delta, \epsilon$ given by the lemma, so
that $\xi$ is admissible in $[0,Q]$ and vanishes at $Q$. We apply
\eqref{E-c} again (making $W_{-} \equiv 0$ which is sufficient for
our estimates). For this purpose, we compute,
\begin{equation*}
(1-2a)(\xi')^2 - 2a \xi \xi'' = \left\{%
\begin{aligned}
& - \frac{a \beta}{R^2} (1 + \frac{t}{R})^{-2\beta -2}, \text{ ~for~ }
t \in [0,R],\\
& ~2 a \alpha (\alpha + 1) \delta^2 \epsilon^2 (1 + \epsilon
Q)^{-\alpha} (1 + \epsilon t)^{-\alpha - 2},
\text{~for~} t \in [R,Q],
\end{aligned}%
\right. 
\end{equation*}
and we obtain,
\begin{equation}\label{E-q3}
\begin{split}
& \frac{a \beta }{R^2} \int_{B(R)} (1 + \frac{r}{R})^{-2 \beta-2} \,
d\mu \le 2 \pi a \big\{ 1 - \sum_{n=1}^{\Nbar}\omega_n (1 +
\frac{t_n}{R})^{-2\beta}
\big\} \\
& + 2 a \alpha (\alpha + 1) \delta^2 \epsilon^2 (1 +
\epsilon Q)^{-\alpha} \int_{C(R,Q)} (1 + \epsilon r)^{-\alpha - 2}
\, d\mu
+ \int_{B(Q)} W_{+} \xi^2(r) \, d\mu,
\end{split}
\end{equation}
where we have chosen $R > t_{\Nbar}$. We fix $R > t_{\Nbar}$, and we
let $Q$ tend to infinity, using the fact that $\delta$ and
$\epsilon$ remain controlled, and that the second term in the
right-hand side of \eqref{E-q3} goes to zero when $Q$ tends to
infinity because $(M,g)$ has $k_a$-subpolynomial volume growth
(Lemmas~\ref{L-B} and \ref{L-N2}). Finally, we obtain,
$$
R^{-2} V(R) \le C \big\{ 1 - \sum_{n=1}^{\Nbar}\omega_n (1 +
\frac{t_n}{R})^{-2\beta} \big\} + \int_{M} W_{+} \, d\mu \le C',
$$
for some constant $C'$ independent of $R$. It follows that $M$ has at most
quadratic volume growth.

\begin{rem}\label{R-6}
When $W_{+} \equiv 0$ and $\chi(M) = 0$ (which corresponds to $\Nbar
= 1$ and $\omega_1=1$), the above estimate gives that $M$ has at
most linear volume growth.
\end{rem}

\subsubsection{Conclusion}

In the three cases (i), (ii) and (iii), we have proved:
\begin{itemize}
    \item[$\diamond$] The surface $M$ has finite topology; when $W_{+} \equiv 0$,
    $M$ is homeomorphic to $\C$ or to $\Cb$ (Step 2).
    \item[$\diamond$] The surface $(M,g)$ has at most quadratic volume growth, and
    hence \cite[Proposition~2.3]{Cas06} is conformally
    equivalent to a closed Riemannian surface with finitely many points removed (Step 3).
    \item[$\diamond$] The function $W$ is integrable and $0 \le 2 \pi
    a \chi(M) + \int_M W \, d\mu$; in particular, when $\chi(M)=0$ and
    $W_{+} \equiv 0$, then $W_{-}\equiv 0$ (Step 2).
    \item[$\diamond$] When $W_{+}\equiv 0$ and $\chi(M)=0$, $(M,g)$ has at most linear
    volume growth .
\end{itemize}

The proof of Assertions~A and B in Theorem~\ref{T-GG} is therefore
complete. \qed

\subsection{Proof of Theorem~\ref{T-GG} -- Assertion C}\label{SS-ptgg2}

By Theorem~\ref{T-GG}, Assertions~A and B, we already know that
$(M,g)$ has finite topology and, more precisely, that $\Nbar$ is
finite. We also know that $W$ is integrable and satisfies the
inequality $0 \le 2\pi a \chi(M) + \int_M W\, d\mu$. Assume that
\begin{equation}\label{E-dag}
0 = 2\pi a \chi(M) + \int_M W\, d\mu.
\end{equation}

In order to prove Assertion~C, we will prove (i) that $(M,g)$ has
subquadratic volume growth, and (ii) that $0 \le a K(x) + W(x)$ for
a.e. $x$ in $M$.

\subsubsection{Proof that $(M,g)$ has subquadratic volume
growth}\label{sss-sqvg}~\\

\boxit{\hbox{Case $a>\4$\,.}} Taking $\xi(t) =
\big(1-\frac{t}{Q}\big)^{\alpha}_{+}$ for some fixed $\alpha$ large
enough. Inequality \eqref{E-c} becomes
$$
\frac{C_{\alpha}}{Q^2} \int_{B(Q)} \big( 1 - \frac{r}{Q}
\big)^{2\alpha-2}\, d\mu \le 2\pi a \big( 1 - \sum_{n=1}^{\Nbar}
\omega_n \xi^2(t_n) \big) + \int_{B(Q)} W \, \xi^2(r) \, d\mu\,.
$$

The term in the right-hand side can be estimated as follows, using
the definition of $\Nbar$ and the assumption \eqref{E-dag},
\begin{equation*}
\begin{split}
2\pi a  & \big( 1 - \sum_{n=1}^{\Nbar} \omega_n \xi^2(t_n) \big) +
\int_{B(Q)} W \, \xi^2(r) \, d\mu \\
& = 2\pi a \big( 1 - \sum_{n=1}^{\Nbar} \omega_n) +
2\pi a \sum_{1}^{\Nbar} \frac{2\alpha t_n \omega_n}{Q} +
o(\frac{1}{Q}) + \int_M W \, \xi^2(r) \, d\mu\\
& \le  \frac{4\pi a\alpha}{Q} \sum_{1}^{\Nbar} t_n
\omega_n + o(\frac{1}{Q})+ 2\pi a \chi(M) + \int_M W \, \xi^2(r) \, d\mu\\
& \le \frac{C'_{\alpha}}{Q} + o(\frac{1}{Q}) + \int_M W
\, \big( \xi^2(r) - 1\big) \, d\mu.
\end{split}
\end{equation*}

It follows that
$$
\frac{C_{\alpha}}{2^{2\alpha-2}} \frac{V(\frac{Q}{2})}{Q^2} \le
\frac{C'_{\alpha}}{Q} + o(\frac{1}{Q}) + \int_M W \, \big( \xi^2(r)
- 1\big) \, d\mu.
$$

Letting $Q$ tend to infinity, and using the dominated convergence
theorem in the right-hand side, we conclude that
$$
\lim_{R\to \infty} \frac{V(R)}{R^2} = 0,
$$
\ie that $(M,g)$ has subquadratic volume growth. \medskip

\boxit{\hbox{Case $a=\4$\,.}} We use the same test function as in
Section~\ref{sss-s3}, Case (ii). We choose $R > t_{\Nbar}$.
Inequality \eqref{E-c} becomes
\begin{equation*}
\begin{split}
& \frac{1}{4R^2} \int_{B(R)}  e^{2(1 - \frac{r}{2R})^2}\, d\mu \le
\frac{\pi}{2}\Big( e^2 - \sum_{n=1}^{\Nbar} \omega_n
e^{2(1-\frac{t_n}{2R})^2}\Big)\\
& + \frac{\alpha^2 \beta^2}{2}e^{-\alpha Q} \int_{C(R,Q)}
e^{-\alpha r}\, d\mu + \int_{B(Q)} W\, \xi^2(r) \, d\mu.
\end{split}
\end{equation*}
Letting $Q$ tend to infinity, using Lemma~\ref{L-N1}, the definition
of $\Nbar$ and the assumption~\eqref{E-dag}, we obtain, for large values of $R$,
\begin{equation*}
\begin{split}
\sqrt{e}\, \frac{V(R)}{4R^2} \le & ~\frac{\pi e^2}{2}\big( 1 -
\sum_{n=1}^{\Nbar} \omega_n e^{-\frac{t_n}{R}(2 -
\frac{t_n}{2R})}\big)\\
& + \int_{B(R)} W e^{2(1-\frac{r}{2R})^2}\, d\mu
+ \beta^2 \int_{M\setminus B(R)} W e^{-2\alpha r}\, d\mu,\\[6pt]
\sqrt{e}\, \frac{V(R)}{4e^2R^2} \le & ~\frac{\pi}{2}\big( 1 -
\sum_{n=1}^{\Nbar} \omega_n \big) + \frac{\pi}{R} \sum_{n=1}^{\Nbar}
t_n \omega_n + o(\frac{1}{R})\\
& + \int_{B(R)} W e^{-\frac{r}{R}(2-\frac{r}{2R})}\,
d\mu + \frac{\beta^2}{e^2} \int_{M\setminus B(R)} W e^{-2\alpha r}\,
d\mu,\\[6pt]
\sqrt{e}\, \frac{V(R)}{4e^2R^2} \le & ~\frac{\pi}{R}
\sum_{n=1}^{\Nbar} t_n \omega_n  + o(\frac{1}{R}) +
\int_{B(R)} W \big( e^{-\frac{r}{R}(2-\frac{r}{2R})} - 1 \big)\\
& + \frac{\beta^2}{e^2} \int_{M\setminus B(R)} W
e^{-2\alpha r}\, d\mu - \int_{M\setminus B(R)} W\, d\mu.\\[6pt]
\end{split}
\end{equation*}

Recall that $\beta$ is uniformy bounded (Lemma~\ref{L-A}) and that
$\alpha > 0$. It follows that $|W|e^{-2\alpha r} \le |W|$, and that $W$ is
integrable (Assertion~A). Using the fact that
$e^{-\frac{r}{R}(2-\frac{r}{2R})} - 1$ tends to zero when $R$ tends
to infinity, and the dominated convergence theorem, we can conclude
that
$$
\lim_{R\to \infty} \frac{V(R)}{R^2} = 0,
$$
\ie that $(M,g)$ has subquadratic volume growth. \medskip

\boxit{\hbox{Case $0 < a <\4$\,.}} We use the same test function as
in Section~\ref{sss-s3}, Case~(iii). We choose $R > t_{\Nbar}$.
Inequality \eqref{E-c} becomes
\begin{equation*}
\begin{split}
& \frac{a \beta }{R^2} \int_{B(R)}  (1 + \frac{r}{R})^{-2\beta -2} \,
d\mu \le 2\pi a \big( 1 - \sum_{n=1}^{\Nbar} \omega_n (1 +
\frac{t_n}{R})^{-2\beta}\big)\\
& + 2a \alpha(\alpha + 1) \delta^2 \epsilon^2 (1+\epsilon
Q)^{-\alpha} \int_{C(R,Q)}(1+\epsilon r)^{-\alpha -2}\, d\mu\\
& + \int_{B(Q)} W \xi^2(r) \, d\mu.
\end{split}
\end{equation*}

Using Lemma~\ref{L-N2}, for $Q$ tending to infinity and for large values of $R$,
we obtain
\begin{equation*}
\begin{split}
\frac{a \beta}{2^{2\beta +2}} \frac{V(R)}{R^2} \le & ~2\pi a \big( 1 -
\sum_{n=1}^{\Nbar} \omega_n\big) + \frac{4\pi a
\beta}{R}\sum_n t_n \omega_n + o(\frac{1}{R}) \\
& + \int_{B(R)} W
(1+\frac{r}{R})^{-2\beta}\, d\mu + \delta^2 \int_{M\setminus B(R)} W
(1+\frac{\epsilon r}{R})^{-2\alpha}\, d\mu\\[6pt]
\frac{a \beta}{2^{2\beta +2}} \frac{V(R)}{R^2} \le & ~\frac{4\pi a
\beta}{R}\sum_n t_n \omega_n  + o(\frac{1}{R}) + \int_{B(R)} W
\big[(1+\frac{r}{R})^{-2\beta} - 1
\big]\, d\mu \\
& - \int_{M\setminus B(R)} W\, d\mu + \delta^2
\int_{M\setminus B(R)} W (1+\frac{\epsilon r}{R})^{-2\alpha}\, d\mu.
\end{split}
\end{equation*}

Recall that $\delta$ is uniformly bounded (Lemma~\ref{L-B}) and that
$|W|(1+\frac{\epsilon r}{R})^{-2\alpha} \le |W|$ which is
integrable. Using the fact that $(1+\frac{r}{R})^{-2\beta} - 1$
tends to zero when $R$ tends to infinity, and the dominated
convergence theorem, we conclude that
$$
\lim_{R\to \infty} \frac{V(R)}{R^2} = 0,
$$
\ie that $(M,g)$ has subquadratic volume growth. \medskip

\subsubsection{Proof that $a K(x) + W(x) \ge 0$ for a.e. $x \in
M$}\label{sss-pos}

Recall that,
\begin{equation}\label{E-pos-a}
0 \le \int_M \Big( |df|^2 + a K f^2 + W f^2\Big) \, d\mu, ~~~~~
\forall f \in \mathrm{Lip}_0(M).
\end{equation}

Fix some $x \in M$, and take the distance function $r$, the geodesic
balls, and the function $\chih(t)$ with respect to this point.
\medskip

\noid According to Assertion~A, $M$ has finite topology and, more
precisely, the function $\chih$ has at most finitely many
discontinuities \ie $\Nbar$ is finite.\medskip

\noid Let $0 < \alpha < 1$ and $0 < R < t_1 < t_{\Nbar} < Q$. Define
the function $\xi$ (with parameters $\alpha, R, Q$) to be

\begin{equation}\label{E-pos-b}
\xi (t) =
\left\{%
\begin{aligned}
& 1 - \alpha \frac{t}{R}, \text{~for~} t \in [0,R],\\
& (1-\alpha) \frac{Q-t}{Q-R}, \text{~for~} t \in [R,Q].\\
\end{aligned}%
\right.
\end{equation}

\noid Use the function $\xi(r)$ to test the positivity condition
\eqref{E-pos-a}. Straightforward computations give
\begin{equation}\label{E-pos-c}
\int_{B(Q)} (\xi')^2(r) \, d\mu = \frac{\alpha^2}{R^2}V(R) + \Big(
\frac{1-\alpha}{Q-R} \Big)^2 \Big( V(Q) - V(R) \Big).
\end{equation}
Applying Lemma~\ref{L-t1} to the ball $B(R)$ and to the set
$C(R,Q)$, another computation yields
\begin{equation*}
\begin{split}
& \int_{B(Q)} K \xi^2(r) \, d\mu \le  - \frac{2\alpha^2}{R^2}V(R) +
2(1-\alpha) \frac{R-\alpha Q}{R(Q-R)} L(R)\\
& + 2\pi \Big( 1 - (1-\alpha)^2 \sum_{n=1}^{\Nbar}\omega_n
(\frac{Q-t_n}{Q-R})^2
\Big)
- 2 (\frac{1-\alpha}{Q-R})^2 \Big( V(Q) - V(R) \Big).
\end{split}%
\end{equation*}

Finally, we obtain that for $a > 0$ and the above choice
\eqref{E-pos-b} of $\xi$,
\begin{equation*}
\begin{split}
& 0 \le \int_{B(Q)} \Big( (\xi')^2(r) + a K \xi^2(r) + W \xi^2(r)
\Big) \, d\mu \le
(1-2a)\alpha^2 \frac{V(R)}{R^2} \\
& + (1-2a) \big( \frac{1-\alpha}{Q-R} \big)^2 \Big( V(Q)
- V(R) \Big) + 2a (1-\alpha) \frac{R-\alpha Q}{Q-R}
\frac{L(R)}{R} \\
& + 2 \pi a \Big( 1 - (1-\alpha)^2 \sum_{n=1}^{\Nbar}
\omega_n (\frac{Q-t_n}{Q-R})^2\Big)\\
&+ \int_{B(R)} W (1 -\alpha \frac{r}{R})^2\, d\mu +
(1-\alpha)^2 \int_{C(R,Q)} W (\frac{Q-r}{Q-R})^2\, d\mu.
\end{split}%
\end{equation*}

\noid The preceding inequality holds for all choices of $\alpha \in
(0,1)$ and $0 < R < t_1  < t_{\Nbar} < Q$.  From
Section~\ref{sss-sqvg}, we know that $(M,g)$ has subquadratic volume
growth. Letting $Q$ tend to infinity, we find that
\begin{equation*}
\begin{split}
0 \le & ~(1-2a)\alpha^2 \frac{V(R)}{R^2} - 2a \alpha (1-\alpha)
\frac{L(R)}{R}
+ 2 \pi a \alpha ( 2 - \alpha )\\
& + 2\pi a (1-\alpha)^2\big(1 -
\sum_{n=1}^{\Nbar}\omega_n \big) + (1-\alpha)^2\int_M W \, d\mu\\
& - (1-\alpha)^2 \int_{B(R)} W\, d\mu + \int_{B(R)} W (1
- \alpha \frac{r}{R})^2\, d\mu,
\end{split}
\end{equation*}
for all $\alpha \in (0,1)$ and $R \in (0, t_1)$.
\medskip

\noid By definition of $\Nbar$ and using \eqref{E-dag}, we have
that,
$$
2\pi a \big(1 - \sum_{n=1}^{\Nbar}\omega_n \big) + \int_M W\, d\mu
\le 0,
$$
and it follows that
\begin{equation*}
\begin{split}
0 \le & ~(1-2a)\alpha^2 \frac{V(R)}{R^2} - 2a \alpha (1-\alpha)
\frac{L(R)}{R} + 2 \pi a \alpha ( 2 - \alpha ) \\
& + \int_{B(R)} W \big[(1 - \alpha \frac{r}{R})^2 -
(1-\alpha)^2\big]\, d\mu.
\end{split}
\end{equation*}
We finally conclude that for all $\alpha \in (0,1)$ and $R \in (0,t_1)$,
\begin{equation}\label{E-pos-f}
\begin{split}
0 \le & ~(1-2a)\alpha \frac{V(R)}{R^2} - 2a (1-\alpha)
\frac{L(R)}{R} + 2 \pi a ( 2 - \alpha )\\
& + \int_{B(R)} W (1 - \frac{r}{R})\big[2-\alpha(1 +
\frac{r}{R})\big]\, d\mu.
\end{split}
\end{equation}

\noid We now use the classical expansions for the length and area of
small geodesic balls $B(x,R)$ when $R$ is small,
\begin{equation}\label{E-pos-g}
\begin{split}
L(R) = 2\pi R \Big( 1 - \frac{K(x)}{6}R^2 + R^2 \epsilon_1(R) \Big) , \\
V(R) = \pi R^2 \Big( 1 - \frac{K(x)}{12}R^2 + R^2 \epsilon_2(R) \Big). \\
\end{split}%
\end{equation}

\noid Plugging \eqref{E-pos-g} into \eqref{E-pos-f}, we find that
\begin{equation}\label{E-pos-h}
\begin{split}
0 \le & ~\pi \alpha + \frac{K(x) R^2}{12} \Big( 8a - (1+6a)
\alpha \Big) \\
& + \pi \, R^2 \Big( (1-2a) \alpha \epsilon_2(R) -
4a(1-\alpha)\epsilon_1(R) \Big)\\
& + \int_{B(R)} W (1 - \frac{r}{R})\big[2-\alpha(1 +
\frac{r}{R})\big]\, d\mu.
\end{split}%
\end{equation}

Letting $\alpha$ tend to zero, we obtain the inequality
$$
0 \le \frac{2\pi}{3}a K(x)R^2 - 4\pi a R^2\epsilon_1(R) +
2\int_{B(R)} W(1-\frac{r}{R})\, d\mu \,,
$$
which holds for all $x \in M$ and $R > 0$ small enough. To finish
the proof, we need to compute the asymptotic expansion of the
integral when $R$ tends to zero. By Lebesgue differentiation theorem
in local coordinates on the surface $M$ \cite[\S 1.7,Theorem 1]{EvGa92},
we have
$$
\int_{B(R)} W\, d\mu = \pi R^2 W(x) + o(R^2)
$$
for a.e. $x\in M$ and $R$ small. Moreover, using normal coordinates $(u,r)\in
S^1\times\R_+$ centered at $x$, the volume form of $M$ reads
$d\mu=\theta_x(u,r)\,du\,dr$, where the density function $\theta_x$
satisfies $\lim_{r\to 0}\frac{\theta_x(u,r)}{r}=1$, and we have
\begin{equation*}
\begin{split}
\frac{1}{R^3}\int_{B(R)} r \, W\, d\mu & =
\frac{1}{R^3}\int_{S^1}\int_0^R W(u,r) \frac{\theta_x(u,r)}{r} r^2\, dr\, du \\
     & = \frac{2}{3}\frac{1}{(R^\frac{3}{2})^2}
        \int_{S^1}\int_0^{R^\frac{3}{2}} W(u,s^\frac{2}{3})
        \frac{\theta_x(u,s^\frac{2}{3})}{s^\frac{2}{3}} s\,ds\, du.
\end{split}
\end{equation*}
Applying Lebesgue differentiation theorem to the function
$W(u,s^\frac{2}{3}) \frac{\theta_x(u,s^\frac{2}{3})}{s^\frac{2}{3}}$
on balls of radii $R^\frac{3}{2}$ we get
$$
\int_{B(R)} r \, W\, d\mu = \frac{2\pi}{3} R^3 W(x) + o(R^3)
$$
for a.e. $x\in M$. Finally, for $R$ small, we have
$$
2\int_{B(R)} W (1-\frac{r}{R})\, d\mu = \frac{2\pi}{3}R^2 W(x) +
o(R^2)
$$
and inequality \eqref{E-pos-h} becomes
$$
0 \le \frac{2\pi}{3}\big( a K(x) + W(x) \big) R^2 + o(R^2),
$$
which concludes the proof that $0 \le a K(x) + W(x)$ for a.e. $x\in
M$. \qed

\subsubsection{Proof that $a K + W \equiv 0$}\label{sss-end}

Recall from Sections~\ref{S-intro} and \ref{SS-ptgg1} that we already deduced
Huber's theorem and Cohn-Vossen's inequality from Theorem~\ref{T-GG},
Assertions A and B.\medskip

The first consequence of the inequality $0 \le a K(x) + W(x)$ for
\emph{a.e.} $x \in M$ is that $a K_{-} \le W_{+}$ a.e. on $M$. Since
$W_{+}$ is integrable by assumption, it follows that $K_{-}$ is
integrable and, by Huber's theorem, that $K$ itself is integrable.
We can then apply Cohn-Vossen's inequality and conclude that $\int_M
K \, d\mu \le 2 \pi \chi (M)$. Finally,
$$
0 \le \int_M (a K + W)\, d\mu \le 2 \pi \chi (M) + \int_M W \, d\mu
= 0.
$$
It follows that $a K + W \equiv 0$ a.e. on $M$. \qed \medskip

\subsection{The volume growth assumptions in Theorem~\ref{T-GG} are
optimal}\label{SS-ptgg3}

To show that the volume growth assumptions in Theorem~\ref{T-GG} are
optimal, we take $W \equiv 0$.\medskip

\noid The example of the hyperbolic plane shows that Assertion~A in
Theorem~\ref{T-G} cannot hold when $a=\4$ without an extra
assumption on $(M,g)$.\medskip

\noid The examples of the hyperbolic planes of curvature $-c^2,
c>0$, show that the assumptions in Case (ii) are optimal. \medskip

\noid Consider the unit disk $\D$ with the conformal metric
$h_{\alpha} = \big( \frac{2}{1-|z|^2} \big)^{2\alpha}|dz|^2$ for
$\alpha \ge 1$. The metric $h_1$ is the hyperbolic metric with
constant curvature $-1$. When $\alpha > 1$, the metric $h_{\alpha}$
is a complete conformal metric on $\D$, with negative curvature. A
simple computation shows that it has polynomial volume growth of
degree $2+\frac{1}{\alpha - 1}$. Given a complete Riemannian surface
$(M,g)$, let $a_{+}(M,g)$ denote the supremum of the numbers $a$
such that $\Delta + a K \ge 0$. This supremum is achieved
\cite[Proposition~1.1]{Cas06} and $a (\alpha):=
a_{+}(\D,h_{\alpha})$ is equal to $\frac{1}{4\alpha}$ by
\cite[Proposition~4.3]{Cas06}. It follows that the volume growth of
$(\D,h_{\alpha})$ is polynomial with degree equal to $2 + \frac{4
a(\alpha)}{1-4a(\alpha)}$. This shows that the volume growth
assumption in Case~(iii) is optimal.

\section{Generalization to finite index operators}\label{S-fio}

Theorems~\ref{T-GG} and \ref{T-G} have their counterparts with the
assumption that the operator $J$ is \emph{non-negative} replaced by
the assumption that the operator $J$ has \emph{finite index}. As a
matter of fact, one can immediately reduce the former case to the
latter by using the following proposition of independent interest.

\begin{prop}\label{P-1}
    Let $(M,g)$ be a complete Riemannian manifold, and let $W$ be a
    locally integrable function on $M$. Then the operator $\Delta + W$
    has finite index if and only if there exists a locally integrable
    function $P$ with compact support such that the operator $\Delta + W
    + P$ is non-negative.
\end{prop}

\pf Assume that $\Delta + W$ has finite index on $C^1_0(M)$. Then
there exists a compact $K \subset M$ such that $\Delta + W$ is
non-negative on $C_0^1(M\setminus K)$. Take $\phi$ to be a smooth
function with compact support, such that $0 \le \phi \le 1$ and
$\phi \equiv 1$ in a compact neighborhood of $K$. Given any $\psi
\in C_0^1(M)$, write $\psi$ as $\psi = \phi \psi + (1-\phi)\psi$. An
easy computation gives,
\begin{equation}\label{E-fip-1}
\begin{split}
& \int_M |d\psi|^2 + W\psi^2 = \int_M |d\big((1-\phi)\psi\big)|^2 +
W\big((1-\phi)\psi\big)^2 \\
& + \int_M W \big( \phi^2 + 2\phi(1-\phi)\big) \psi^2
- \2 \int_M \psi^2 \Delta \big((1-\phi)^2\big) - \int_M
\psi^2 |d\phi|^2\\
& + 2 \int_M \phi (1-\2\phi)|d\psi|^2\,.\\
\end{split}
\end{equation}
Because $\Delta + W$ is non-negative in $M\setminus K$, and because
of our choice of $\phi$, the first and fourth terms in the
right-hand side of \eqref{E-fip-1} are non-negative. The other terms
can be written as $-\int_M P \psi^2$, where the function $P$ is
defined by
\begin{equation}\label{E-fip-2}
P :=  |d\phi|^2 - \Delta \big( \phi (1 - \2 \phi)\big)
 - W\phi^2 - 2 \phi (1-\phi)W.\\
\end{equation}

Recall that $W$ is locally integrable and that $\phi$ is smooth with
compact support. It follows that $P$ is locally integrable, with
compact support. By \eqref{E-fip-1}, the operator $\Delta + W + P$
is non-negative on $C_0^1(M)$, as stated.\medskip

\noid Assume that there exists a function $P$, which is locally
integrable with compact support, such that $\Delta + W + P$ is
non-negative on $C_0^1(M)$. Let $K$ be a compact neighborhood of the
support of $P$. Then,
$$
0 \le \int_M |d\psi|^2 + W\psi^2 + P\psi^2 = \int_M |d\psi|^2 +
W\psi^2,
$$
for any $\psi \in C_0^1(M\setminus K)$, and this means that $\Delta
+ W$ is non-negative on $C_0^1(M\setminus K)$. By a result of
B.~Devyver \cite{Dev10}, this implies that $\Delta + W$ has finite
index on $C_0^1(M)$. \qed

\begin{thm}\label{T-GFI}
Let $(M,g)$ be a complete non-compact Riemannian surface. Let $W$ be
a locally integrable function on $M$. Assume that the function $W_+$
is integrable on $M$, and that the operator $\Delta + a K +W$ has
finite index. Assume furthermore that either,
\begin{itemize}
\item[(i)] $a > \4$, or
\item[(ii)] $a = \4$, and $(M,g)$ has subexponential volume growth, or
\item[(iii)] $a \in (0,\4)$, and $(M,g)$ has $k_a$-subpolynomial volume
growth, with $k_a = 2 + \frac{4a}{1-4a}$.
\end{itemize}
Then $W$ is integrable and $(M,g)$ is conformally equivalent to a
closed Riemannian surface with finitely many points removed.
\end{thm}

\pf By Proposition~\ref{P-1}, there exists some function $P$, with
compact support, such that $\Delta + aK + W + P$ is non-negative. As
$(W+P)_+$ is still integrable, we can apply Theorem \ref{T-GG} to
conclude that $(M,g)$ is conformally equivalent to a closed
Riemannian surface with finitely many points removed and that $W+P$
is integrable. Since $P$ has compact support, we also have that $W$
is integrable. \qed

\section{Other proofs}\label{S-proofs}

\subsection{Proofs of Theorem~\ref{T-G} and \ref{T-C}}\label{SS-ptgc}

To prove Theorem~\ref{T-G} and \ref{T-C}, it suffices to follow the
proof of Theorem~\ref{T-GG}, Assertions~A and B, making $W_{+}
\equiv 0$, $q = W_{-}$, and to use the Remarks~\ref{R-1}-\ref{R-4}.

\subsection{Proof of Theorem~\ref{T-R}}\label{SS-ptr}

Let us first consider the case of the sphere with constant curvature
$\alpha^2$, $M_0 = S^2(\alpha^2)$. In the sequel, the subscript $0$
refers to $M_0$. Let $J_0 = \Delta_0 + a \alpha^2 - c$. The operator
$J_0$ is non-negative in the ball $B_0(R)$ if and only if the first
Dirichlet eigenvalue of the Laplacian $\Delta_0$ in this ball
satisfies $\lambda_1\big(B_0(R)\big) \ge c - a\alpha^2$. Since $c
\ge (a+2)\alpha^2$, it follows that $J_0$ non-negative in the ball
$B_0(R)$ implies that $\lambda_1\big( B_0(R)\big) \ge \lambda_1\big(
B_0(\frac{\pi}{2\alpha})\big)$, and hence that $R \le
\frac{\pi}{2\alpha}$, because $\lambda_1\big( B_0(R)\big)$ is a
decreasing function of $R$. If $J_0 \ge 0$ in
$B_0(\frac{\pi}{2\alpha})$, then $c=(a+2)\alpha^2$, since all
previous inequalities become equalities. Recall that the first
Dirichlet eigenfunction for the Laplacian $\Delta_0$ in the
hemisphere $B_0(\frac{\pi}{2\alpha})$ is $\cos (\alpha r_0)$, up to
a scaling factor, where $r_0$ is the distance function to a point on
the sphere. \medskip

\emph{Proof of Theorem~\ref{T-R}}. Recall that this theorem is of a
local nature. We first state a lemma.

\begin{lem}\label{L-TR}
Let $(M,g)$ be a Riemannian surface. Assume that the curvature
satisfies $K \le \alpha^2$ for some $\alpha > 0$. Let $J$ be the
operator $J = \Delta + a K - c$, with $a \in [\2, \infty)$ and $c
\ge (a + 2) \alpha^2$. Assume furthermore that the ball $B(x,
\frac{\pi}{2\alpha})$ is contained in $M$, for some $x \in M$. Then
the least eigenvalue of the operator $J$ with Dirichlet boundary
conditions in this ball is non-positive. If $J$ has least Dirichlet
eigenvalue $0$ in the ball $B(x, \frac{\pi}{2\alpha})$, then
$c=(a+2)\alpha^2$, $K \equiv \alpha^2$ and $B(x,
\frac{\pi}{2\alpha})$ is covered by the hemisphere
$S_{+}^2(\alpha^2)$ in the sphere with constant curvature
$\alpha^2$.
\end{lem}

Clearly, the lemma implies the theorem. Indeed, Assertion~(i)
follows from the lemma and from the monotonicity of eigenvalues with
respect to domain inclusion. Assertion~(ii) follows
immediately.\medskip

\emph{Proof of the lemma}.

\noid First observe that we can reduce to the case $a \in [\2,2]$.
Indeed, if $a
> 2$, then for any $a' \in [\2,2]$, we can write
$$
\Delta + aK-c = \Delta + a'K + (a-a')K - c \le \Delta + a'K -c',
$$
where $c' = c + (a'-a)\alpha^2 \ge (a'+2)\alpha^2$. Moreover, if
$c'=(a'+2)\alpha^2$, then $c=(a+2)\alpha^2$.\medskip

\noid Assume that $a \in [\2,2]$. Let $A := \frac{\pi}{2\alpha}$.
Because $K\le \alpha^2$, the map $\exp_x : T_xM \to M$ is a local
diffeomorphism on the ball $D(0,A)$. Let $\tilde{g} = \exp_x^{*}g$
be the pulled-back metric to $T_xM$. Let $\mu_1$ be the least
Dirichlet eigenvalue of $\Delta + a K - c$ in $B(x,A)$. Then,
$\Delta + a K - c - \mu_1 \ge 0$ in $B(x,A)$ and hence, there exists
a positive function $u : B(x,A) \to \R$ such that $(\Delta + a K - c
- \mu_1)u = 0$, see \cite{FCSc80}. Let $\tilde{u} = u\circ \exp_x$.
Because $\exp_x$ is a local isometry, we have $(\tilde{\Delta} + a
\tilde{K} - c - \mu_1)\tilde{u} = 0$ and hence the least Dirichlet
eigenvalue $\tilde{\mu}_1$ of $\tilde{\Delta} + a \tilde{K} - c$
satisfies $\tilde{\mu}_1 \ge \mu_1$. To show that $\mu_1$ is
non-positive, it suffice to show that $\tilde{\mu}_1$ is
non-positive. We have reduced to the simply-connected case. \medskip

\noid We now work in the simply-connected disk $D(0,A)$, with a
metric (also denoted) $g$ such that $K \le \alpha^2$. We denote by
$L(r)$ the length of $\partial D(0,r)$ for this metric and we let
$L_0(r)$ be the corresponding length on the sphere, $L_0 (r) = 2\pi
\frac{\sin (\alpha r)}{\alpha}$. By Bishop's comparison theorem, we
have that
\begin{equation}\label{E-ma}
L (r) \ge L_0(r).
\end{equation}

We now use Pogorelov's trick. Let $\xi : [0,A] \to \R$ be a $C^2$
function such that $\xi(0)=1$ and $\xi(A)=0$. We compute the
quadratic form $\mathcal{Q}$ associated with $J = \Delta + a K - c$ on the
function $\xi(r)$, where $r$ is the geodesic distance to $0$ in
$D(0,A)$. We also introduce the total curvature of $D(0,r)$,
\begin{equation}\label{E-mb}
G(r) = \int_{D(r)} K \, d\mu ,
\end{equation}
with respect to the Riemannian measure in $D(0,A)$. Applying the
co-area formula,
\begin{equation*}
\begin{split}
\mathcal{Q}(\xi(r)) & = \int_{D(0,A)} \big( |d\xi(r)|^2 + (aK-c) \xi^2(r)
\big) \, d\mu \\
& = \int_{0}^A \big( (\xi')^2 - c \xi^2\big)(t) L(t) \, dt + a \int_0^A G'(t)
\xi^2(t) \, dt,
\end{split}
\end{equation*}
and we compute the second integral in the right-hand side by
integration by parts,
\begin{equation}\label{E-mc}
\mathcal{Q}(\xi(r)) = \int_{0}^A \Big( (\xi')^2 - c \xi^2\Big)(t) L(t) \, dt - a
\int_0^A G(t) (\xi^2)'(t) \, dt.
\end{equation}

By the Gauss-Bonnet formula, we have $G(t) = 2 \pi - L'(t)$ so that
\eqref{E-mc} becomes, after another integration by parts,
$$
\mathcal{Q}(\xi(r)) = \int_{0}^A \Big( (\xi')^2 - c \xi^2 - a (\xi^2)''\Big)(t) L(t)
\, dt + 2\pi a.
$$

Finally, we obtain
\begin{equation}\label{E-md}
\mathcal{Q}(\xi(r)) = \int_{0}^A \Big( (1-2a) (\xi')^2 - 2 a \xi \xi'' - c
\xi^2 \Big)(t) L(t) \, dt + 2\pi a,
\end{equation}
for any function $\xi : [0,A] \to \R$ which is $C^2$ and such that
$\xi(0)=1$ and $\xi(A)=0$. \medskip

We now use the test function $\eta(r) = \cos (\alpha r)$ in formula
\eqref{E-md}, where $r$ is the Riemannian distance to the center of
the ball, \ie we transplant the first eigenfunction of the
hemisphere to a function on the ball $D(0,A)$.
$$
\mathcal{Q}(\eta(r)) = \int_0^A \Big( (1-2a) \alpha^2 \sin^2(\alpha t) +
(2a\alpha^2 - c) \cos^2(\alpha t) \Big) L(t) \, dt + 2\pi a,
$$
and hence
\begin{equation}\label{E-me}
\mathcal{Q}(\eta(r)) \le \int_0^A \Big( (1-2a) \alpha^2 \sin^2(\alpha t) +
(a-2)\alpha^2 \cos^2(\alpha t) \Big) L(t) \, dt + 2\pi a,
\end{equation}
where we have used the fact that $c \ge (2+a)\alpha^2$. Recall that
$a \in [\2,2]$. Using the inequality \eqref{E-ma}, we find that
\begin{equation}\label{E-mf}
\mathcal{Q}(\eta(r)) \le 2\pi a
+ \int_0^A \big( (1-2a) \alpha^2 \sin^2(\alpha t) +
(a-2)\alpha^2 \cos^2(\alpha t) \big) L_0(t) \, dt.
\end{equation}

The right-hand side of \eqref{E-mf} is zero because this is the
value of the quadratic form of the operator $J_0 = \Delta -
2\alpha^2$ on the hemisphere $S^2_+(\alpha^2) =
B_0(\frac{\pi}{2\alpha})$. We conclude that $\mathcal{Q}(\eta(r)) \le 0$ and
hence that the least Dirichlet eigenvalue of $\Delta + a K - c$ in
$D(0,A)$ is non-positive, as stated in the lemma. If this eigenvalue
is zero, then $\mathcal{Q}(\eta(r)) = 0$, and we must have equality in both
\eqref{E-me} and \eqref{E-mf}, \ie $c = (2+a)\alpha^2$ and $L(t)
\equiv 2 \pi \frac{\sin (\alpha t)}{\alpha}$. We then deduce that
$G(t) \equiv 2\pi (1 - \cos (\alpha t))$. Since $K \le \alpha^2$,
integrating $K$ we find that $K \equiv \alpha^2$ and hence we
conclude that $D(0,A) = S^2_{+}(\alpha^2)$. This proves the lemma.
\qed \medskip

As a corollary of Theorem~\ref{T-R}, we obtain Mazet's estimates.

\begin{cor}\label{C-R}
Let $(M,g) \looparrowright (\Mh, \gh)$ be an isometric immersion
with constant mean curvature $H$ in a simply connected space form
with constant sectional curvature $\kappa$. Assume furthermore that
$H^2+\kappa > 0$ and that the immersion is (strongly) stable. Then,
$$
d_g(x, \partial M) \le \frac{\pi}{2 \sqrt{H^2+\kappa}},
$$
where $d_g(x, \partial M)$ is the distance from $x \in M$ with
respect to the metric $g$ to the boundary of $M$, with equality if
and only if $M$ is the hemisphere of a sphere of mean curvature $H$
in $\Mh$.
\end{cor}

\pf The Jacobi operator of the immersion is $J = \Delta - |A|^2 -
\widehat{\mathrm{Ric}}(n)$, where $A$ is the second fundamental form
of the immersion and $n$ the unit normal along the immersion. By the
Gauss equation, we find that $J = \Delta + 2K - 4(H^2 + \kappa)$ and
that $K = H^2 + \kappa - \4 (k_1-k_2)^2$, where $k_i$ are the
principal curvatures. We can apply Theorem~\ref{T-R} with $a=2$ and
$\alpha^2 = H^2 + \kappa$. For the equality case, note that equality
implies that $M$ is totally umbilic. \qed \medskip

\textbf{Remarks}.
\begin{enumerate}
\item This corollary provides a unified proof of Theorem~3.1 and Corollary~3.2 in
Mazet's paper \cite{Maz09}, without using Lawson's correspondence.
\item The proof of Theorem~\ref{T-R} is simpler than that of
\cite[Theorem~3.1]{Maz09}, but it uses the same idea which
goes back to A.~Pogorelov \cite{Pog81}.
\item It is not clear whether the assumptions $a\ge \frac{1}{2}$
and $c \ge (a+2)\alpha^2$ are sharp.
\end{enumerate}

\newcommand{\RMIauthor}{\newblock}
\newcommand{\RMIpaper}{\newblock}
\newcommand{\RMIjournal}{\newblock}
\newcommand{\RMIbook}{\newblock}

\vspace*{15mm}

\begin{footnotesize}
\noindent\begin{minipage}{0.4\textwidth}
    \begin{flushleft}
    {\normalsize Pierre B\'{e}rard}\\
    Universit\'{e} de Grenoble\\
    Institut Fourier (\textsc{ujf-cnrs})\\
    B.P. 74\\
    38402 Saint Martin d'H\`{e}res Cedex\\
    France\\
    \verb+pierrehberard@gmail.com+\\
    \end{flushleft}
\end{minipage}
\hfill
\begin{minipage}{0.45\textwidth}
    \begin{flushleft}
    {\normalsize Philippe Castillon}\\
    Universit\'{e} Montpellier II\\
    D\'{e}pt des sciences math\'{e}matiques CC 51\\
    I3M (\textsc{umr 5149})\\
    34095 Montpellier Cedex 5\\
    France\\
    \verb+philippe.castillon@univ-montp2.fr+\\
    \end{flushleft}
\end{minipage}
\end{footnotesize}


\begin{thebibliography}{99}
\parskip=0.7pt


\bibitem{BerCas12a}
\RMIauthor{B\'{e}rard, P. and Castillon, P.}
\RMIpaper{Remarks on J.~Espinar\lowercase{'s} ``Finite index operators
on surfaces''}
\RMIjournal{arXiv:1204.1604v1}

\bibitem{BerCas12b}
\RMIauthor{B\'{e}rard, P. and Castillon, P.}
\RMIpaper{Spectral positivity and Riemannian coverings}
\RMIjournal{Bull. London Math. Soc.} \textbf{45} (2013), no. 5, 1041--1048

\bibitem{CaPe79}
\RMIauthor{Carmo do, M. and Peng, C-K.}
\RMIpaper{Stable complete minimal surfaces in $\R^3$ are planes}
\RMIjournal{Bull. Amer. Math. Soc.} \textbf{1} (1979), no. 6, 903--906


\bibitem{Cas06}
\RMIauthor{Castillon, P.}
\RMIpaper{An inverse spectral problem on surfaces}
\RMIjournal{Comment. Math. Helv.} \textbf{81} (2006), 271--286

\bibitem{CoMi02}
\RMIauthor{Colding, T. and Minicozzi, W.}
\RMIpaper{Estimates for parametric elliptic integrands}
\RMIjournal{Internat. Math. Res. Notices} \textbf{6} (2002), 291--297

\bibitem{Dev10}
\RMIauthor{Devyver, B.}
\RMIpaper{On the finiteness of the Morse index for Schr\"{o}dinger operators}
\RMIjournal{Manuscripta Math.} \textbf{139} (2012), no. 1-2, 249--271

\bibitem{Esp09}
\RMIauthor{Espinar, J.M.}
\RMIpaper{Finite index operators on surfaces}
\RMIjournal{J. Geometric Analysis} \textbf{23} (2013), no. 1, 415--437

\bibitem{Esp11}
\RMIauthor{Espinar, J.M.}
\RMIpaper{Rigidity of stable cylinders in three-manifolds}
\RMIjournal{Proc. Amer. Math. Soc.} \textbf{140} (2012), no. 5, 1769--1775

\bibitem{EsRo10}
\RMIauthor{Espinar, J. and Rosenberg, H.}
\RMIpaper{A Colding-Minicozzi inequality and its applications}
\RMIjournal{Trans. Amer. Math. Soc.} \textbf{363} (2011), no. 5, 2447--2465

\bibitem{EvGa92}
\RMIauthor{Evans, L.C. and Gariepy, R.F.}
\RMIbook{Measure theory and fine properties of functions}
CRC Press, Boca Raton, FL, 1992

\bibitem{FCSc80}
\RMIauthor{Fischer-Colbrie, D. and Schoen, R.}
\RMIpaper{The structure of complete stable minimal surfaces in
$3$-manifolds of non-negative scalar curvature}
\RMIjournal{Comm. Pure Applied Math.} \textbf{33} (1980), 199--211

\bibitem{Fia40}
\RMIauthor{Fiala, F.}
\RMIpaper{Le probl\`{e}me des isop\'{e}rim\`{e}tres sur les surfaces ouvertes \`{a}
courbure positive}
\RMIjournal{Comment. Math. Helv.} \textbf{13} (1940-1941), 293--346

\bibitem{GuLa86}
\RMIauthor{Gulliver, R. and Lawson, H.B.}
\RMIpaper{The structure of minimal hypersurfaces near a singularity}
\RMIjournal{Proc. Symp. Pure Math.} \textbf{44} (1986), 213--237

\bibitem{Har64}
\RMIauthor{Hartman, P.}
\RMIpaper{Geodesic parallel coordinates in the large}
\RMIjournal{Amer. J. Math.} \textbf{86} (1964), no. 4, 705--727

\bibitem{Kaw88}
\RMIauthor{Kawai, S.}
\RMIpaper{Operator {$\Delta - aK$} on surfaces}
\RMIjournal{Hokkaido Math. J.} \textbf{17} (1988), 147--150

\bibitem{Maz09}
\RMIauthor{Mazet, L.}
\RMIpaper{Optimal length estimates for stable CMC surfaces in
$3$-space forms}
\RMIjournal{Proc. Amer. Math. Soc.} \textbf{137} (2009), no. 8, 2761--2765

\bibitem{MePeRo06}
\RMIauthor{Meeks, W.H., P\'{e}rez, J. and Ros, A.}
\RMIpaper{Liouville-type properties for embedded minimal surfaces}
\RMIjournal{Comm. Analysis and Geometry} \textbf{4} (2006), no. 4, 703--723

\bibitem{MePeRo08}
\RMIauthor{Meeks, W.H., P\'{e}rez, J. and Ros, A.}
\RMIpaper{Stable constant mean curvature surfaces}
In \RMIbook{Handbook of Geometric Analysis}, Vol.~1. Lizhen Ji, Peter Li,
Richard Schoen and Leon Simon Editors. International Press, 2008

\bibitem{Miy93}
\RMIauthor{Miyaoka, R.}
\RMIpaper{$L^2$ harmonic $1$-forms on a complete stable minimal hypersurface}
In \RMIbook{Geometry and global analysis} Ed. T.~Kotake, S.~Nishikawa,
R.~Schoen. Sendai 1993, pp. 289--293

\bibitem{Pog81}
\RMIauthor{Pogorelov, A.V.}
\RMIpaper{On the stability of minimal surfaces}
\RMIjournal{Soviet Math. Dokl.} \textbf{24} (1981), 274--276

\bibitem{Rei10}
\RMIauthor{Reiris, M.}
\RMIpaper{Geometric relations of stable minimal surfaces and
applications}
\texttt{arXiv:1002.3274v1}

\bibitem{ScYa82}
\RMIauthor{Schoen, R. and Yau, S.T.}
\RMIpaper{Complete three dimensional manifolds with positive Ricci curvature
and scalar curvature}
In \RMIbook{Seminar on differential geometry}
Ed. S.T.~Yau. Annals of Math. Studies {\bf 102} (1982), 209-228

\bibitem{ShTa89}
\RMIauthor{Shiohama, K. and Tanaka, M.}
\RMIpaper{An isoperimetric problem for infinitely connected complete
open surfaces}
In \RMIbook{Geometry of Manifolds (Mastumoto, 1988)},
Prospect. Math. {\bf 8}. Academic Press 1989, 317-343

\bibitem{ShTa93}
\RMIauthor{Shiohama, K. and Tanaka, M.}
\RMIpaper{The length function of geodesic parallel circles}
In \RMIbook{Progress in differential geometry}, (K. Shiohama Ed.). Adv. Stud. Pure
Math. {\bf 22}, Math. Soc. Japan, Tokyo 1993, 299-308





\end{thebibliography}
\end{document}